% OneTape.tex 
% joint paper with Dan Seabold
% to be submitted to the Mathematical Logic Quarterly
% July, 1999
\font\fifteenrm=cmr10 scaled\magstep2 % this is all really just 14.4pt
\font\fifteeni=cmmi10 scaled\magstep2
\font\fifteensy=cmsy10 scaled\magstep2
\font\fifteenbf=cmbx10 scaled\magstep2
\font\fifteentt=cmtt10 scaled\magstep2
\font\fifteenit=cmti10 scaled\magstep2
\font\fifteensl=cmsl10 scaled\magstep2
\font\fifteenam=msam10 scaled\magstep2
\font\fifteenbm=msbm10 scaled\magstep2
\font\fifteenex=cmex10 scaled\magstep2
\font\fifteensc=cmcsc10 scaled\magstep2 
\font\twelverm=cmr10 at 12pt
\font\twelvei=cmmi10 at 12pt
\font\twelvesy=cmsy10 at 12pt
\font\twelvebf=cmbx10 at 12pt
\font\twelvett=cmtt10 at 12pt
\font\twelveit=cmti10 at 12pt
\font\twelvesl=cmsl10 at 12pt
\font\twelveam=msam10 at 12pt
\font\twelvebm=msbm10 at 12pt
\font\twelveex=cmex10 at 12pt
\font\twelvesc=cmcsc10 at 12pt
\font\elevenrm=cmr10 scaled\magstephalf % this is really 10.95 pt
\font\eleveni=cmmi10 scaled\magstephalf
\font\elevensy=cmsy10 scaled\magstephalf
\font\elevenbf=cmbx10 scaled\magstephalf
\font\eleventt=cmtt10 scaled\magstephalf
\font\elevenit=cmti10 scaled\magstephalf
\font\elevensl=cmsl10 scaled\magstephalf
\font\elevenam=msam10 scaled\magstephalf
\font\elevenbm=msbm10 scaled\magstephalf
\font\elevenex=cmex10 scaled\magstephalf
\font\elevensc=cmcsc10 scaled\magstephalf
\font\tenrm=cmr10
\font\teni=cmmi10
\font\tensy=cmsy10
\font\tenbf=cmbx10
\font\tentt=cmtt10
\font\tenit=cmti10
\font\tensl=cmsl10
\font\tenam=msam10
\font\tenbm=msbm10
\font\tenex=cmex10
\font\tensc=cmcsc10
\font\ninerm=cmr9
\font\ninei=cmmi9
\font\ninesy=cmsy9
\font\ninebf=cmbx9
\font\ninett=cmtt9
\font\nineit=cmti9
\font\ninesl=cmsl9
\font\nineam=msam9
\font\ninebm=msbm9
\font\nineex=cmex9
\font\ninesc=cmcsc9
\font\eightrm=cmr8
\font\eighti=cmmi8
\font\eightsy=cmsy8
\font\eightbf=cmbx8
\font\eighttt=cmtt8
\font\eightit=cmti8
\font\eightsl=cmsl8
\font\eightam=msam8
\font\eightbm=msbm8
\font\eightex=cmex8
\font\eightsc=cmcsc8
\font\sevenrm=cmr7
\font\seveni=cmmi7
\font\sevensy=cmsy7
\font\sevenbf=cmbx7

\font\sevenam=msam7
\font\sevenbm=msbm7

\font\sixrm=cmr6
\font\sixi=cmmi6
\font\sixsy=cmsy6

\font\sixit=cmti7 at 6pt

\font\sixam=msam6
\font\sixbm=msbm6

\font\fiverm=cmr5
\font\fivei=cmmi5
\font\fivesy=cmsy5

\font\fiveam=msam5
\font\fivebm=msbm5

\font\fourrm=cmr5 at 4pt
\font\fouri=cmmi5 at 4pt
\font\foursy=cmsy5 at 4pt

\font\fouram=msam5 at 4pt
\font\fourbm=msbm5 at 4pt

\skewchar\twelvei='177 \skewchar\eleveni='177\skewchar\teni='177
\skewchar\ninei='177 \skewchar\eighti='177\skewchar\seveni='177 
\skewchar\sixi='177 \skewchar\fivei='177 \skewchar\fouri='177
\skewchar\twelvesy='60 \skewchar\elevensy='60 \skewchar\tensy='60
\skewchar\ninesy='60 \skewchar\eightsy='60 \skewchar\sevensy='60 
\skewchar\sixsy='60 \skewchar\fivesy='60 \skewchar\foursy='60
\newfam\itfam
\newfam\slfam
\newfam\bffam
\newfam\ttfam
\newfam\scfam
\newfam\amfam
\newfam\bmfam
\def\eightbig#1{{\hbox{$\left#1\vbox to 6.5pt{}\voidright $}}}
\def\eightBig#1{{\hbox{$\left#1\vbox to 7.5pt{}\voidright $}}}
\def\eightbigg#1{{\hbox{$\left#1\vbox to 10pt{}\voidright $}}}
\def\eightBigg#1{{\hbox{$\left#1\vbox to 13pt{}\voidright $}}}
\def\ninebig#1{{\hbox{$\left#1\vbox to 7.5pt{}\voidright $}}}
\def\nineBig#1{{\hbox{$\left#1\vbox to 8.5pt{}\voidright $}}}
\def\ninebigg#1{{\hbox{$\left#1\vbox to 11.5pt{}\voidright $}}}
\def\nineBigg#1{{\hbox{$\left#1\vbox to 14.5pt{}\voidright $}}}
\def\tenbig#1{{\hbox{$\left#1\vbox to 8.5pt{}\voidright $}}}
\def\tenBig#1{{\hbox{$\left#1\vbox to 9.5pt{}\voidright $}}}
\def\tenbigg#1{{\hbox{$\left#1\vbox to 12.5pt{}\voidright $}}}
\def\tenBigg#1{{\hbox{$\left#1\vbox to 16pt{}\voidright $}}}
\def\elevenbig#1{{\hbox{$\left#1\vbox to 9pt{}\voidright $}}}
\def\elevenBig#1{{\hbox{$\left#1\vbox to 10.5pt{}\voidright $}}}
\def\elevenbigg#1{{\hbox{$\left#1\vbox to 14pt{}\voidright $}}}
\def\elevenBigg#1{{\hbox{$\left#1\vbox to 17.5pt{}\voidright $}}}
\def\twelvebig#1{{\hbox{$\left#1\vbox to 10pt{}\voidright $}}}
\def\twelveBig#1{{\hbox{$\left#1\vbox to 11pt{}\voidright $}}}
\def\twelvebigg#1{{\hbox{$\left#1\vbox to 15pt{}\voidright $}}}
\def\twelveBigg#1{{\hbox{$\left#1\vbox to 19pt{}\voidright $}}}
\def\fifteenbig#1{{\hbox{$\left#1\vbox to 12pt{}\voidright $}}}
\def\fifteenBig#1{{\hbox{$\left#1\vbox to 13.5pt{}\voidright $}}}
\def\fifteenbigg#1{{\hbox{$\left#1\vbox to 18pt{}\voidright $}}}
\def\fifteenBigg#1{{\hbox{$\left#1\vbox to 23pt{}\voidright $}}}
\def\voidright{\right.\nulldelimiterspace=0pt \mathsurround=0pt }
\def\fifteenpoint{
  \textfont0=\fifteenrm \scriptfont0=\twelverm \scriptscriptfont0=\tenrm
  \def\rm{\fam0 \fifteenrm}%
  \textfont1=\fifteeni \scriptfont1=\twelvei \scriptscriptfont1=\teni
  \textfont2=\fifteensy \scriptfont2=\twelvesy \scriptscriptfont2=\tensy
  \textfont3=\fifteenex \scriptfont3=\fifteenex \scriptscriptfont3=\fifteenex
  \def\it{\fam\itfam\fifteenit}\textfont\itfam=\fifteenit
  \def\sl{\fam\slfam\fifteensl}\textfont\slfam=\fifteensl
  \def\bf{\fam\bffam\fifteenbf}\textfont\bffam=\fifteenbf 
    \scriptfont\bffam=\twelvebf\scriptscriptfont\bffam=\tenbf
  \def\tt{\fam\ttfam\fifteentt}\textfont\ttfam=\fifteentt
  \def\sc{\fam\scfam\fifteensc}\textfont\scfam=\fifteensc
  \def\am{\fam\amfam\fifteenam}\textfont\amfam=\fifteenam
    \scriptfont\amfam=\twelveam\scriptscriptfont\amfam=\tenam
  \def\bm{\fam\bmfam\fifteenbm}\textfont\bmfam=\fifteenbm
    \scriptfont\bmfam=\twelvebm\scriptscriptfont\bmfam=\tenbm
  \baselineskip=21pt \rm
  \let\big=\fifteenbig\let\Big=\fifteenBig\let\bigg=\fifteenbigg
  \let\Bigg=\fifteenBigg}
\def\twelvepoint{
  \textfont0=\twelverm \scriptfont0=\ninerm \scriptscriptfont0=\sevenrm
  \def\rm{\fam0 \twelverm}%
  \textfont1=\twelvei \scriptfont1=\ninei \scriptscriptfont1=\seveni
  \textfont2=\twelvesy \scriptfont2=\ninesy \scriptscriptfont2=\sevensy
  \textfont3=\twelveex \scriptfont3=\twelveex \scriptscriptfont3=\twelveex
  \def\it{\fam\itfam\twelveit}\textfont\itfam=\twelveit
  \def\sl{\fam\slfam\twelvesl}\textfont\slfam=\twelvesl
  \def\bf{\fam\bffam\twelvebf}\textfont\bffam=\twelvebf 
    \scriptfont\bffam=\ninebf\scriptscriptfont\bffam=\sevenbf
  \def\tt{\fam\ttfam\twelvett}\textfont\ttfam=\twelvett
  \def\sc{\fam\scfam\twelvesc}\textfont\scfam=\twelvesc
  \def\am{\fam\amfam\twelveam}\textfont\amfam=\twelveam
    \scriptfont\amfam=\nineam\scriptscriptfont\amfam=\sevenam
  \def\bm{\fam\bmfam\twelvebm}\textfont\bmfam=\twelvebm
    \scriptfont\bmfam=\ninebm\scriptscriptfont\bmfam=\sevenbm
  \baselineskip=17.8pt \rm 
  \def\looselineskip{\baselineskip=18.5pt plus 1.8pt}%
  \def\tightlineskip{\baselineskip=16.5pt}%
  \def\verytightlineskip{\baselineskip=15pt}%
  \let\big=\twelvebig\let\Big=\twelveBig\let\bigg=\twelvebigg
  \let\Bigg=\twelveBigg  }
\def\elevenpoint{
  \textfont0=\elevenrm \scriptfont0=\ninerm \scriptscriptfont0=\sixrm
  \def\rm{\fam0 \elevenrm}%
  \textfont1=\eleveni \scriptfont1=\ninei \scriptscriptfont1=\sixi
  \textfont2=\elevensy \scriptfont2=\ninesy \scriptfont2=\sixsy 
  \textfont3=\elevenex \scriptfont3=\elevenex \scriptfont3=\elevenex
  \def\it{\fam\itfam\elevenit}\textfont\itfam=\elevenit
  \def\sl{\fam\slfam\elevensl}\textfont\slfam=\elevensl
  \def\bf{\fam\bffam\elevenbf}\textfont\bffam=\elevenbf
  \def\tt{\fam\ttfam\eleventt}\textfont\ttfam=\eleventt
  \def\sc{\fam\scfam\elevensc}\textfont\scfam=\elevensc
  \def\am{\fam\amfam\elevenam}\textfont\amfam=\elevenam
    \scriptfont\amfam=\nineam\scriptscriptfont\amfam=\sixam
  \def\bm{\fam\bmfam\elevenbm}\textfont\bmfam=\elevenbm
    \scriptfont\bmfam=\ninebm\scriptscriptfont\bmfam=\sixbm
  \baselineskip=15.1pt \rm
  \def\looselineskip{\baselineskip=16pt plus 1.5pt}%
  \def\tightlineskip{\baselineskip=14pt}%
  \def\verytightlineskip{\baselineskip=13pt}%
  \let\big=\elevenbig\let\Big=\elevenBig\let\bigg=\elevenbigg
  \let\Bigg=\elevenBigg  }
\def\tenpoint{
  \textfont0=\tenrm \scriptfont0=\eightrm \scriptscriptfont0=\fiverm
  \def\rm{\fam0 \tenrm}%
  \textfont1=\teni \scriptfont1=\eighti \scriptscriptfont1=\fivei
  \textfont2=\tensy \scriptfont2=\eightsy \scriptfont2=\fivesy 
  \textfont3=\tenex \scriptfont3=\tenex \scriptfont3=\tenex
  \def\it{\fam\itfam\tenit}\textfont\itfam=\tenit
  \def\sl{\fam\slfam\tensl}\textfont\slfam=\tensl
  \def\bf{\fam\bffam\tenbf}\textfont\bffam=\tenbf
  \def\tt{\fam\ttfam\tentt}\textfont\ttfam=\tentt
  \def\sc{\fam\scfam\tensc}\textfont\scfam=\tensc
  \def\am{\fam\amfam\tenam}\textfont\amfam=\tenam
    \scriptfont\amfam=\eightam \scriptscriptfont\amfam=\fiveam
  \def\bm{\fam\bmfam\tenbm}\textfont\bmfam=\tenbm
    \scriptfont\bmfam=\eightbm \scriptscriptfont\bmfam=\fivebm
  \baselineskip=14pt\rm
  \def\looselineskip{\baselineskip=14.8pt plus1.5pt}
  \def\tightlineskip{\baselineskip=12.6pt}%
  \def\verytightlineskip{\baselineskip=13pt}%
  \let\big=\tenbig\let\Big=\tenBig\let\bigg=\tenbigg\let\Bigg=\tenBigg  }
\def\ninepoint{
  \textfont0=\ninerm \scriptfont0=\sevenrm \scriptscriptfont0=\fourrm
  \def\rm{\fam0 \ninerm}%
  \textfont1=\ninei \scriptfont1=\seveni \scriptscriptfont1=\fouri
  \textfont2=\ninesy \scriptfont2=\sevensy \scriptfont2=\foursy 
  \textfont3=\nineex \scriptfont3=\nineex \scriptfont3=\nineex
  \def\it{\fam\itfam\nineit}\textfont\itfam=\nineit
  \def\sl{\fam\slfam\ninesl}\textfont\slfam=\ninesl
  \def\bf{\fam\bffam\ninebf}\textfont\bffam=\ninebf
  \def\tt{\fam\ttfam\ninett}\textfont\ttfam=\ninett
  \def\sc{\fam\scfam\ninesc}\textfont\scfam=\ninesc
  \def\am{\fam\amfam\nineam}\textfont\amfam=\nineam
    \scriptfont\amfam=\nineam\scriptscriptfont\amfam=\fouram
  \def\bm{\fam\bmfam\ninebm}\textfont\bmfam=\ninebm
    \scriptfont\bmfam=\ninebm\scriptscriptfont\bmfam=\fourbm
  \baselineskip=12.6pt\rm
  \def\tightlineskip{\baselineskip=11.5pt}
  \let\big=\ninebig\let\Big=\nineBig\let\bigg=\ninebigg
  \let\Bigg=\nineBigg  }
\def\eightpoint{
  \textfont0=\eightrm \scriptfont0=\fiverm \scriptscriptfont0=\fourrm
  \def\rm{\fam0 \eightrm}%
  \textfont1=\eighti \scriptfont1=\fivei \scriptscriptfont1=\fouri
  \textfont2=\eightsy \scriptfont2=\fivesy \scriptfont2=\foursy 
  \textfont3=\eightex \scriptfont3=\eightex \scriptfont3=\eightex
  \def\it{\fam\itfam\eightit}\textfont\itfam=\eightit
  \def\sl{\fam\slfam\eightsl}\textfont\slfam=\eightsl
  \def\bf{\fam\bffam\eightbf}\textfont\bffam=\eightbf
  \def\tt{\fam\ttfam\eighttt}\textfont\ttfam=\eighttt
  \def\sc{\fam\scfam\eightsc}\textfont\scfam=\eightsc
  \def\am{\fam\amfam\eightam}\textfont\amfam=\eightam
    \scriptfont\amfam=\eightam\scriptscriptfont\amfam=\fouram
  \def\bm{\fam\bmfam\eightbm}\textfont\bmfam=\eightbm
    \scriptfont\bmfam=\eightbm\scriptscriptfont\bmfam=\fourbm
  \baselineskip=11.2pt \rm
  \def\tightlineskip{\baselineskip=10.4pt}
  \let\big=\eightbig\let\Big=\eightBig\let\bigg=\eightbigg
  \let\Bigg=\eightBigg  }

\twelvepoint
\nopagenumbers
\hsize=6in\vsize=8.8in

\parskip=1pt plus 1pt

\newif\ifSpecialhead\Specialheadfalse
\newbox\specialheadbox

\def\specialhead #1\par{\Specialheadtrue\setbox\specialheadbox=\hbox{#1}}
\headline={{\ifSpecialhead\box\specialheadbox\global\Specialheadfalse\else
     \ifnum\pageno<0{\hfill\quad{\twelvebf\folio}}%
     \else\ifnum\pageno<2\hfill
     \else\hfill\twelvepoint\sc\firstmark\quad{\twelvebf\folio}\fi\fi\fi}}

\def\title#1\par{\medskip{\def\cr{\hfil\par\hfil}\hfil\fifteenbf#1\hfil\par}\medskip}
\def\subtitle#1\par{\centerline{\fifteenrm #1}\medskip}
\def\author#1\par{\medskip{\def\cr{\hfil\par\hfil\twelvesc}\fifteensc\hfil#1\hfil\par}}
\def\authors#1\par{\medskip\fifteensc\center#1\par}
\def\center#1\par{{\def\cr{\hfil\par\hfil}\hfil#1\hfil\par}}
\def\abstract.#1\par{\message{Abstract.}%
                    \medskip{\narrower\narrower\tenpoint\tightlineskip
                        \noindent{\bf Abstract.}#1\par}\medskip\noindent}
\def\tinyabstract.#1\par{\message{Abstract.}%
                    \medskip{\narrower\narrower\eightpoint\tightlineskip
                        \noindent{\bf Abstract.}#1\par}\medskip\noindent}
\def\bigabstract.#1\par{\message{Abstract.}%
                         \medskip{\narrower\narrower\tightlineskip
                         \noindent{\bf Abstract. }#1\par}\medskip\noindent}
\def\acknowledgement#1\par{\footnote{}{#1}}
\def\sectionskip{\Goodbreak\vskip 25pt plus 15pt minus 5pt}
\def\secnumber{\ifquiet
               \else\ifNoSections
                    \else\sectionsymbol\the\secno\quad\fi\fi}
\def\section#1\par{ \NoSectionsfalse\par\sectionskip\proofdepth=0\claimno=0
 \ifquiet\else\advance\secno by1\fi\toks0={#1}
 \immediate\write16{\ifquiet\else Section \the\secno\space\fi
                    \the\toks0}%
 \mark{\secnumber #1}%
 {\fifteenpoint\bf\noindent\secnumber #1}\nobreak\bigskip\quietoff
 \nobreak\noindent}
\def\quiet{\quiettrue}

\def\quietoff{\ifQUIET\else\quietfalse\fi}
\newif\ifquiet
\newif\ifQUIET
\newif\ifNoSections
\newcount\claimtype
\newcount\secno
\newcount\claimno
\newcount\subclaimno
\newcount\subsubclaimno
\newcount\subsubsubclaimno
\newcount\proofdepth
\def\subclaimnumber{\ifquiet\else\ifcase\subclaimno\or A\or B\or C\or D\or E\or
     F\or G\or H\or I\or J\or K\or L\or M\or N\or O\or P\fi\fi}
\def\subsubclaimnumber{\ifquiet\else\ifcase\subsubclaimno\or i\or ii\or iii\or 
   iv\or v\or vi\or vii\or viii\or ix\or x\or xi\or xii\or xiii\or xiv\fi\fi}
\def\subsubsubclaimnumber{\ifquiet\else\ifcase\subsubsubclaimno\or a\or b\or 
   c\or d\or e\or f\or g\or \or h\or i\or j\or k\or l\or m\or n\or o\fi\fi}
\def\claimtag{\ifquiet\else
  \ifNoSections
    \ifcase\proofdepth\the\claimno%
    \or\the\claimno.\subclaimnumber
    \or\the\claimno.\subclaimnumber.\subsubclaimnumber
    \or\the\claimno.\subclaimnumber.\subsubclaimnumber
                                                .\subsubsubclaimnumber\fi
  \else
    \ifcase\proofdepth\the\secno.\the\claimno
    \or\the\secno.\the\claimno.\subclaimnumber
    \or\the\secno.\the\claimno.\subclaimnumber.\subsubclaimnumber
    \or\the\secno.\the\claimno.\subclaimnumber.\subsubclaimnumber
                                                .\subsubsubclaimnumber\fi\fi\fi}
\secno=0\claimno=0\proofdepth=0\subclaimno=0\subsubclaimno=0\subsubsubclaimno=0
\NoSectionstrue
\newbox\qedbox
\def\claimname{\ifcase\claimtype Theorem\or Lemma\or Claim\or Corollary\or
               Question\or Definition\or Remark\or Conjecture\fi}
\def\preclaimskip{\removelastskip
    \ifcase\claimtype\goodbreak\vskip 8pt plus 4pt minus 2pt
                  \or\goodbreak\vskip 6pt plus 4pt minus 1pt
                  \or\goodbreak\vskip 5pt plus 4pt minus 1pt
                  \or\goodbreak\vskip 8pt plus 4pt minus 2pt
                  \or\vskip 7pt plus 4pt minus 2pt
                  \or\vskip 7pt plus 4pt minus 2pt
                  \or\vskip 7pt plus 4pt minus 2pt
                  \or\goodbreak\vskip 8pt plus 4pt minus 2pt\fi}
\def\postclaimskip{\ifcase\claimtype         \vskip 4pt plus 2pt minus 2pt
                                          \or\vskip 3pt plus 2pt minus 2pt
                                          \or\vskip 2pt plus 2pt minus 1pt
                                          \or\vskip 4pt plus 2pt minus 2pt
                                          \or\vskip 1pt plus 2pt 
                                          \or\vskip 4pt plus 4pt 
                                          \or\vskip 3pt plus 2pt
                                          \or\vskip 4pt plus 2pt minus 2pt\fi}
\def\claimfont{\ifcase\claimtype
                  \sl\or\sl\or\sl\or\sl\or\sl\or\rm\or\rm\or\sl\fi}
\def\advancetag{\ifcase\proofdepth\advance\claimno by1
                               \or\advance\subclaimno by1
                               \or\advance\subsubclaimno by1
                               \or\advance\subsubsubclaimno by1\fi}
\def\sayclaim#1.#2 #3\par{\ifquiet\else\advancetag\fi
    \preclaimskip\setbox1=\hbox{#1}\setbox2=\hbox{#2}%
    \toks0={#1 }
    \immediate\write16{\ifdim\wd1>0pt\the\toks0
                       \else\claimname\space\fi \claimtag.}%
    \vbox{\noindent
    {\bf\ifdim\wd1=0pt \claimname\else #1\fi\ifquiet.\else\ \claimtag{\ifNoSections.\fi}\fi}%
    \enspace{\ifdim\wd2>0pt\sc #2\enspace\fi}%
    {\claimfont #3\par}}\postclaimskip\quietoff}
\def\theorem{\claimtype=0\sayclaim}
\def\lemma{\claimtype=1\sayclaim}

\def\corollary{\claimtype=3\sayclaim}
\def\question{\claimtype=4\sayclaim}

\def\point#1. #2\par{\item{\rm #1.}#2\par}
\def\points#1\cr\par{\medskip\vbox{\let\cr=\point\point#1\par}\par}
\def\df{\it}
\def\prooffont{}
\def\proofsize{}%\ifcase\proofdepth\or\elevenpoint\or\tenpoint\or\ninepoint\fi}
\def\proofindent{}%\ifcase\proofdepth\or\or\narrower\or\narrower\fi}
\def\proofskip{\badbreak\ifcase\claimtype    \vskip 3pt plus 2pt minus 2pt
                                          \or\vskip 2pt plus 2pt minus 2pt
                                          \or\vskip 1pt plus 2pt minus 1pt
                                          \or\vskip 3pt plus 2pt minus 2pt
                                          \or\vskip 1pt plus 2pt 
                                          \or\vskip 2pt plus 4pt 
                                          \or\vskip 1pt plus 2pt
                                          \or\vskip 3pt plus 2pt minus 2pt\fi}

\def\Goodbreak{\vskip0pt plus.5in\penalty-1000\vskip0pt plus-.5in}
\def\goodbreak{\penalty-500}
\def\badbreak{\penalty500}
\def\Badbreak{\penalty1000}
\def\proof{\message{proof}\removelastskip\Badbreak\proofskip\begingroup
  \advance\proofdepth by1
  \setbox\qedbox=\hbox{\halmos\raise2pt\hbox{\fiverm\claimname}}%
  \prooffont\proofsize\proofindent\noindent{\bf Proof: }}
\def\proofof#1:{\message{proof}\removelastskip\Badbreak\proofskip\begingroup
  \advance\proofdepth by1
  \setbox\qedbox=\hbox{\halmos\raise2pt\hbox{\fiverm#1}}%
  \prooffont\proofsize\proofindent\noindent{\bf Proof of #1: }}
\def\cite[#1]{[{\tenrm{#1}}]\message{[#1]}}
\edef\ref#1{\expandafter\global\expandafter\edef#1{\noexpand\claimtag}}
\newwrite\notes
\openout\notes=\jobname.notes
% for the following macros, see page 377 of the TeXbook
\long\def\unexpandedwrite#1#2{\def\finwrite{\write#1}%
   {\aftergroup\finwrite\aftergroup{\sanitize#2\endsanity}}}
\def\sanitize{\futurelet\next\sanswitch}
\let\stoken=\space
\def\sanswitch{\ifx\next\endsanity
  \else\ifcat\noexpand\next\stoken\aftergroup\space\let\next=\eat
   \else\ifcat\noexpand\next\bgroup\aftergroup{\let\next=\eat
    \else\ifcat\noexpand\next\egroup\aftergroup}\let\next=\eat
     \else\let\next=\copytoken\fi\fi\fi\fi \next}
\def\eat{\afterassignment\sanitize \let\next= }
\long\def\copytoken#1{\ifcat\noexpand#1\relax\aftergroup\noexpand
  \else\ifcat\noexpand#1\noexpand~\aftergroup\noexpand\fi\fi
  \aftergroup#1\sanitize}
\def\endsanity\endsanity{}

\def\note#1#2{\hbox to2in{\strut#1\quad\dotfill\quad#2}}
\def\boxit#1{\setbox4=\hbox{\kern1pt#1\kern1pt}
  \hbox{\vrule\vbox{\hrule\kern1pt\box4\kern1pt\hrule}\vrule}}
\def\halmos{\hbox{\am\char'3}} 
\def\qed#1\par{\message{.                                }\setbox1=\hbox{#1}%
  \ifdim\wd1>0pt\setbox\qedbox=\hbox{\halmos\raise2pt\hbox{\fiverm #1}}\fi
  \kern5pt\lower 2pt\hbox{\box\qedbox}\proofskip\goodbreak\endgroup}

\def\sectionsymbol{\S}

\def\a{\alpha}
\def\b{\beta}
\def\d{\delta}
\def\s{\sigma}

\def\I1{\mathop{\hbox{\sc i}_1}}
\def\w{\omega}

\def\Q{{\mathchoice{\hbox{\bm Q}}{\hbox{\bm Q}}
         {\hbox{\tenbm Q}}{\hbox{\sevenbm Q}}}}
\def\R{{\mathchoice{\hbox{\bm R}}{\hbox{\bm R}}
         {\hbox{\tenbm R}}{\hbox{\sevenbm R}}}}

\def\N{{\mathchoice{\hbox{\bm N}}{\hbox{\bm N}}
         {\hbox{\tenbm N}}{\hbox{\sevenbm N}}}}

\def\unifto{\buildrel\lower 7pt\hbox{$\to$}\over\to}

\def\from{\mathbin{\vbox{\baselineskip=3pt\lineskiplimit=0pt
                         \hbox{.}\hbox{.}\hbox{.}}}}

\def\in{\mathrel{\mathchoice{\raise 
1pt\hbox{$\scriptstyle\cal\char'62$}}
         {\raise 1pt\hbox{$\scriptstyle\cal\char'62$}}
         {\raise .5pt\hbox{$\scriptscriptstyle\cal\char'62$}}
         {\hbox{$\scriptscriptstyle\cal\char'62$}}}\penalty700{}}
\def\ni{\mathrel{\mathchoice{\raise 1pt\hbox{$\scriptstyle\cal\char'63$}}
                   {\raise 1pt\hbox{$\scriptstyle\cal\char'63$}}
                   {\raise .5pt\hbox{$\scriptscriptstyle\cal\char'63$}}
                   {\hbox{$\scriptscriptstyle\cal\char'63$}}}\penalty700}
\def\of{\mathrel{\mathchoice{\raise 1pt\hbox{$\scriptstyle\subseteq$}}
                   {\raise 1pt\hbox{$\scriptstyle\subseteq$}}
                   {\raise .5pt\hbox{$\scriptscriptstyle\subseteq$}}
                   {\hbox{$\scriptscriptstyle\subseteq$}}}}
\def\fo{\mathrel{\mathchoice{\raise 1pt\hbox{$\scriptstyle\supseteq$}}
                   {\raise 1pt\hbox{$\scriptstyle\supseteq$}}
                   {\raise .5pt\hbox{$\scriptscriptstyle\supseteq$}}
                   {\hbox{$\scriptscriptstyle\supseteq$}}}}
\def\notin{\mathrel{\mathchoice
  {\raise 1pt\hbox{\rlap{$\scriptstyle\;|$}$\scriptstyle\cal\char'62$}}
  {\raise 1pt\hbox{\rlap{$\scriptstyle\kern2pt 
          |$}$\scriptstyle\cal\char'62$}}
  {\raise .5pt\hbox{\rlap{$\scriptscriptstyle\, |$}$\scriptscriptstyle
      \cal\char'62$}}
  {\hbox{\rlap{$\scriptscriptstyle\, |$}$\scriptscriptstyle
     \cal\char'62$}}}%
  \penalty700}

\def\and{\mathrel{\kern1pt\&\kern1pt}}

\def\cross{\times}

\def\lte{\mathrel{\scriptstyle\leq}}

\def\[#1]{\left[\vphantom{\bigm|}#1\right]}
\def\<#1>{\langle\,#1\,\rangle}

\def\concat{\mathbin{{\,\hat{ }\,}}}

\def\st{\mid}
\def\seq<#1>{{\def\st{\mid\penalty650}\left<\,#1\,\right>}}

\def\set#1{\{\,{#1}\,\}}

\def\th{{\hbox{\fiverm th}}}

\def\lttheta{{\raise 1pt\hbox{$\scriptstyle<$}\theta}}

\def\I1{\mathop{\hbox{\sc i}_1}}

\def\ilte{\lte_{\scriptscriptstyle\infty}}
\def\iequiv{\equiv_{\scriptscriptstyle\infty}}
\def\concat{\mathop{*}}

\title Infinite time Turing machines with only one tape

\author Joel David Hamkins\cr
            The City University of New York\cr
	{\tentt http://www.math.csi.cuny.edu/$\sim$hamkins}\cr

\author Daniel Evan Seabold\cr
            Hofstra University\cr
            {\tentt http://www.hofstra.edu/$\sim$matdes}\cr            

\acknowledgement The research of the first author has been supported in
part by a grant from the PSC-CUNY Research Foundation. Keywords: one-tape 
infinite time Turing machines, supertask computation. Subject Codes: 0D310, 0D360.

\abstract. Infinite time Turing machines with only one tape are in
many respects fully as powerful as their multi-tape cousins. In
particular, the two models of machine give rise to the same class of
decidable sets, the same degree structure and, at least for functions
$f\from\R\to\N$, the same class of computable
functions. Nevertheless, there are infinite time computable functions
$f:\R\to\R$ that are not one-tape computable, and so the two models 
of infinitary computation 
are not equivalent. Surprisingly, the class of one-tape computable
functions is not closed under composition; but closing it under
composition yields the full class of all infinite time computable
functions. Finally, every ordinal which is clockable by an infinite
time Turing machine is clockable by a one-tape machine, except certain
isolated ordinals that end gaps in the clockable ordinals. 

Infinite time Turing machines, introduced in \cite[HamLew$\infty$a],
extend the usual operation of Turing machines into transfinite ordinal
time. By doing so, they provide a model for supertask computations,
computations involving infinitely many steps, and set the stage for a
mathematical analysis of what is possible in principle to achieve via
supertasks. For example, it is easy to see that every arithmetic set is
decidable by such machines; a bit more sophistication reveals that every
$\Pi^1_1$ set and more is supertask decidable. A rich degree structure has
emerged on the class of reals and sets of reals, stratified by two natural
jump operators. For this and more analysis we refer the reader to the
small but rapidly growing body of literature on the subject:
\cite[HamLew$\infty$a], \cite[HamLew$\infty$b], \cite[Wel$\infty$a],
\cite[Wel$\infty$b] and \cite[Wel98].  

Let us review how the machines work. 
Using a three-tape Turing machine model, with separate
input, scratch and output tapes, an infinite time Turing machine
progresses through the successor stages of computation just as an
ordinary Turing machine does, according to the rigid instructions of a
finite program running with finitely many states. The new behavior appears
at a limit stage: the head is reset to the initial starting position; the
machine is placed in the special {\it limit} state, just another of the
finitely many states; and the values in the cells of the tape are updated
by computing the $\limsup$ of the previous cell values. With the limit
stage configuration thus completely specified, the machine simply
continues computing. If the {\it halt} state is eventually reached, the
machine gives as output whatever is written on the output tape. Since
there is plenty of time for the machines to handle infinite binary input
and output, the natural context for the machines is therefore Cantor Space
$2^\w$, which we denote by $\R$ and refer to as the set of reals. Thus,
the machines provide supertask notions of computability for partial
functions $f\from\R\to\R$ as well as notions of decidability and
semi-decidability for sets of reals $A\of\R$. 

For convenience, the machines are defined with three tapes.  This 
mechanical configuration allows one to keep the input separate from
the scratch work and the output. For example, many arguments in
\cite[HamLew$\infty$a] and \cite[HamLew$\infty$b] begin by regarding
the scratch tape as an infinite list of infinite tapes, kept track of by
means of G\"odel pair coding.

\def\cell#1{\boxit{\hbox to 17pt{\strut\hfil$#1$\hfil}}}
\def\boxit#1{\setbox4=\hbox{\kern2pt#1\kern2pt}\hbox{\vrule\vbox{\hrule\kern2pt\box4\kern2pt\hrule}\vrule}}
\def\col#1#2#3{\hbox{\vbox{\baselineskip=0pt\parskip=0pt\cell#1\cell#2\cell#3}}}
$$\halign{\qquad\hfill#\hfill&\quad\quad\quad\hfill#\hfill\cr
\hbox{\raise 5pt\hbox{\vbox{\hbox to .6in{\tenit\hfill input: }\hbox to.6in{\tenit\hfill scratch: }\hbox to.6in{\tenit\hfill output: }}}\col110\col001\col110\lower2pt\vbox{\hbox{\strut\ninerm\hskip 3pt Head}\boxit{\col011}}\col101\col010\raise 5pt\hbox{\vbox{\hbox{\ $\cdots$}\hbox{\ $\cdots$}\hbox{\ $\cdots$}}}}&\raise20pt\hbox{\cell1\cell0\cell1\lower2pt\vbox{\hbox{\strut\ninerm\hskip3pt Head}\boxit{\cell1}}\cell0\raise5pt\hbox{\ $\cdots$}}\cr
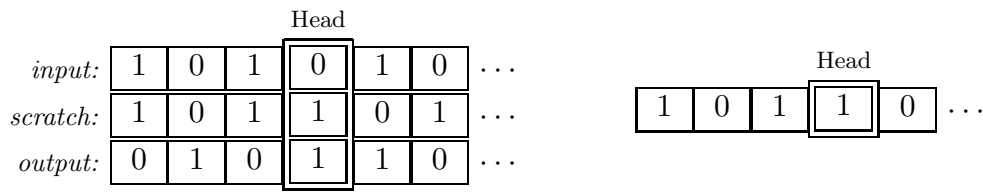
{\eightrm \phantom{scratch: }An ordinary infinite time Turing machine}&{\eightrm A one-tape infinite time Turing machine}\cr}$$

It seems natural to wonder, {\it what if one uses a machine with only one tape?}  When such a machine
begins, the one tape is filled completely with the input, and the
machine embarks on its computation according to the usual rules for
infinite time Turing machines. If, after some transfinite
amount of time, the machine attains its halt state, then whatever is
written on the one tape is the output of the machine. One naturally
obtains the notions of a {\df one-tape computable} function, a {\df
one-tape decidable} set, a {\df one-tape clockable ordinal} and so on, in
analogy with the original theory. The fundamental question we would 
like to consider in this paper is:

\quiet\question Primary Question. Are the infinite time Turing machines with only one tape fully as powerful as the 
original three-tape variety?

The answer is delicately mixed. On the one hand, we show that the two
kinds of machines give rise to exactly the same class of decidable
sets, the same degree structure and, at least for functions whose
range is contained in $\set{0,1}$, $\N$, $\Q$ or $\set{1}\cross\R$,
the same class of computable functions. Thus, the one-tape machines
seem fully as powerful as their three-tape cousins. On the other hand,
we show that there are computable functions $f:\R\to\R$ which are not
computed by any one-tape machine; indeed, the class of one-tape
computable functions is not even closed under composition. So the
one-tape machines are less powerful. Nevertheless, every computable
function is in a precise sense nearly computable by a 
one-tape machine, and the closure of the class of one-tape computable
functions under composition yields the full class of all infinite time
computable functions. Here are the main theorems: 

\quiet\theorem Main Theorems. 
\points 1. A set is decidable if and only if it is one-tape decidable.\cr
	2. A function $f\from\R\to\R$ whose range is not dense in $\R$ is computable if and
only if it is one-tape computable.\cr 
	3. There is a computable function $f:\R\to\R$ which is not one-tape computable.\cr
	4. The class of one-tape computable functions is not closed under composition; closing it under composition yields the class of all computable functions.\cr
	5. Every clockable ordinal is one-tape clockable, except certain isolated ordinals that end gaps in the clockable ordinals.\cr

We will follow the notation and terminology of \cite[HamLew$\infty$a]
and \cite[HamLew$\infty$b]. In particular, by such unadorned terms as
{\df computable} and {\df decidable} we mean computable and decidable
by three-tape infinite time Turing machines. Sometimes,
for emphasis, we will use the term {\df three-tape computable} and so
on. Since we denote by $\R$ the Cantor space $2^\omega$, let us denote
by $\N$ the set of sequences of the form $\<11\cdots 1100\cdots>$,
that is, those with an initial block of $1$s and then all $0$s, and by
$\Q$ the set of sequences that are eventually $0$. If $a$ is a real,
let $a_i$ be the $i^\th$ digit of $a$, so that
$a=\<a_0a_1a_2\cdots>$. By $1\concat a$ we mean the real
$\<1a_0a_1a_2\cdots>$, with a $1$ concatenated to the front of $a$. If
$f$ is a function with range contained in $\R$, then by $1\concat f$
we mean the function $x\mapsto 1\concat f(x)$. We write $f\from A\to
B$ to mean that $f$ is a partial function from $A$ to $B$.  
\section The one-tape machines seem fully powerful

In this section we prove a variety of theorems which build towards the
conclusion that the one-tape infinite time Turing machines are fully
as powerful as their three-tape cousins. We begin by proving that for
a large class of functions, the two models yield the same notion of
computability. Since a three-tape machine can easily simulate one-tape
machine computations, the difficult direction of this argument will be
to simulate a three-tape computation on a one-tape machine. Our basic
strategy for doing so can be divided into three broad phases:

$$\hbox{\tenit Stretch Input}\quad\rightarrow\quad
	\hbox{\tenit Simulate Computation}\quad\rightarrow\quad
	\hbox{\tenit Compress Output}$$
$$\hbox{\eightrm The Three Phases of a Simulated Computation}$$

We describe the middle phase first, since it explains the need for
the other two phases. We view the one tape of a one-tape Turing
machine as divided into blocks of three cells, each block representing
one column on a three-tape machine; that is, each block represents 
one cell each from the input, scratch and output tapes of
a simulated three-tape machine.  While the head of a three-tape machine
can read and write on all three of these cells at once, our simulation
will take up to nine steps to accomplish the equivalent
effect. Specifically, given a three-tape program $p$, there is a
one-tape program $q$ which will simulate the operation of $p$ by
reading the three cells of the current block (in three steps), then
writing on those three cells in the way that $p$ would on the three
tapes (in another three steps), and then finally (in a final three
steps) moving to the next block to the right or the left, accordingly
as $p$ moves left or right in that situation. Thus, by this means each step of
computation of $p$ on a three-tape machine is simulated by $q$ in 
nine steps of one-tape computation. Actually, by combining the steps between
the three modes, a moment's thought shows that in fact seven 
steps suffice. Note that limit stages are simulated
correctly by this scheme because the three tapes are directly represented
cell-for-cell on the one tape, and so the $\limsup$ operation is the same
in both cases. So we have proved the following:  

\lemma Simulation Lemma. For every infinite time Turing machine
program $p$ there is a one-tape program $q$ which simulates $p$ in the
sense that if $p$ halts on input $a$ with $x$, $y$ and $z$,
respectively, on the input, scratch and output tapes, then $q$ halts
on input $\<a_000a_100a_200\cdots>$ with output
$\<x_0y_0z_0x_1y_1z_1\cdots>$. Each step of the $p$ computation
corresponds to seven steps in the $q$ computation. 
\ref\Simulate 

In order to fully simulate a computation by means of the Simulation
Lemma, however, one must transform the input $a$ for $p$ into the
input $\<a_000a_100\cdots>$ for $q$. This explains the need for the
first phase in the simulation strategy mentioned above. Specifically,
the first phase of a full simulation consists of computing the {\df
stretch} function $s:a\mapsto \<a_000a_100a_200\cdots>$, in order to
set up the input configuration correctly for the Simulation Lemma.  

\lemma Stretch Lemma. The stretch function is one-tape computable.
\ref\Stretch

\proof Our procedure takes $\omega^{2}+1$ many steps. The basic idea
is that with each limit we stretch the input by two additional cells,
so that after infinitely many limits we have stretched it fully.  
From the {\it start} state, the machine writes $0$ in the second cell and
$1$ in the third cell and, by means of states, remembers the values of
the two cells it has overwritten. Next, in $\omega$ many steps, it
moves every digit of the input after $a_0$ to the right two cells.
That is, it writes $a_1$ in the fourth cell, $a_2$ in the fifth cell,
and so on.  The $1$ in the 
third cell serves as a movable marker to remind the machine after each
limit stage just how far the input has been stretched.  The second cell
is used to set a flag that will tell us when the task has been
completed.

At stage $\omega$, the machine therefore finds itself in the {\it
limit} state with tape reading $\< a_001a_1a_2a_3\cdots >$. The
machine now ``flashes'' the second cell by writing a $1$ and then a
$0$, shifts every digit starting with $a_{2}$ two 
cells to the right, and moves the marker from the third to the sixth
cell.  Thus at stage $\omega+\omega$ the tape reads $\<a_000a_{1}01a_{2}a_{3}a_{4}\cdots>$, and at 
$\omega+\omega+\omega$, the machine appears as follows:
$$\lower2pt\vbox{\hbox{\strut\hskip5pt\sixit limit}\boxit{\cell{a_0}}}\vtop{\baselineskip=0pt\cell0\hbox to22pt{\hfil\sixrm flash\hfil}\hbox to21pt{\hfil\sixrm cell\hfil}}\cell0\cell{a_1}\cell0\cell0\cell{a_2}\cell0\vtop{\baselineskip=0pt\cell1\rlap{\hskip-1.5pt\sixrm movable}\rlap{\hskip-2.5pt\sixrm\ marker}}\cell{a_3}\cell{a_4}\cell{a_5}\raise5pt\hbox{\ $\cdots$}$$

The machine continues in this fashion: at each limit stage, it
flashes the second cell, moves the marker, and stretches the input an
additional two cells.  After $\omega^{2}$ steps, the input has been
stretched to the desired form, and the machine for the
first time sees a $1$ in the second cell at a limit state. The machine
then erases that cell and  halts, giving the output
$\<a_000a_100a_200\cdots>$, as desired.  \qed 

Supposing we had instead defined the stretch function on $a$ to give the output
$\<00a_000a_100a_2\cdots>$,  then the procedure above would take only $\omega^{2}$
steps, since the flashing flag in this case could be placed on the first cell. In our 
simulations, however, we prefer to preserve the order of the tapes, and so we have 
defined the stretch function in order to put the input on the first cell of each block. 

The final, third phase is the one that compresses the simulated output
from the representation $\<x_0y_0z_0x_1y_1z_1\cdots>$ to just
$z=\<z_0z_1z_2\cdots>$, the contents of the simulated output tape. Let
us call the function which accomplishes this the {\df compression}
function. It is easy to see that the compression function is one-tape
computable in $\omega$ many steps, so let us record that fact here: 

\lemma Compression Lemma. The compression function is one-tape computable.
\ref\Compression

Thus, we have shown that each of the three phases of our strategy for
simulating a three-tape computation with a one-tape machine is
possible by itself. We may now put these three steps together to
obtain the following theorem.  

\theorem Overly Hopeful Theorem. If the class of one-tape computable
functions is closed under composition, then it is the same as the
class of all infinite time computable functions. 
\ref\OverlyHopeful

\proof Given any infinite time computable function $f$, let $g$ be the
simulation function as in the Simulation Lemma \Simulate, and let $s$
and $c$ be the stretch and compression functions, respectively. Since
we have proved that each of these functions is one-tape computable,
the result follows from the simple observation that the original
function can be expressed as the triple composition $f=c\circ g\circ
s$, corresponding to the three phases of our strategy.\qed

Why have we called this theorem overly hopeful?  It is because of the 
simple fact, a fact we were very 
suprised to discover, that both the hypothesis and the conclusion of the theorem
are false. The class 
of one-tape computable functions is not closed under composition, and it
is not the same as the class of all infinite time computable
functions. It is not possible in the general case to
combine the three phases of our strategy and simulate an arbitrary
computable function $f\from\R\to\R$. 
 
We can, nevertheless, use our strategy to simulate the computation of
functions in many cases. So let us begin with a general connection
between computable and one-tape computable functions. Recall that
$1\concat f$ is the function which concatenates a $1$ to the front of
every value of $f$. 

\theorem. A function $f$ is computable if and only if $1\concat f$ is
one-tape computable.
\ref\Almost

Before proving the theorem we would like to pause in order to 
address a small matter not addressed in 
the original definition of infinite time Turing machines in
\cite[HamLew$\infty$a]. The question is, namely, what happens when an infinite 
time Turing machines attempts to move left from the left-most cell? In practice, this question may be avoided in the three-tape context, because one may easily 
put a flag or some such information marking the left-most cell on one of the tapes  
in such a way that 
any program can be replaced with a {\it tidy} program, one which never attempts to 
move left from the left-most cell. But when we are simulating these machines on one-tape machines, we want to simulate arbitrary programs, and not only 
the tidy ones, so it seems best to 
adopt a specific convention. Here, therefore, we adopt the convention that if, after having read the left-most cell, written on that cell and changed to a new state, the machine is then directed by the program to move left, then in fact the head simply stays on the left-most 
cell without moving at all. Thus, attempting to move left from the left-most cell is the same as not moving at all, except the machine does not necessarily realize that it has not actually moved left. The reason that we bring up this issue here is that when we simulate the
operation of a three-tape program $p$ with a one-tape machine we will need to simulate this situation appropriately. Our simulation set-up must therefore be able to know which blocks of cells are the left-most simulated cells. In the Simulation Lemma \Simulate, these cells are the left-most cells on the one-tape machine, and so no additional care needs to be taken there. But often, we will have 
additional flags in front of the simulated blocks of cells, and so at these times we will need to take additional care. We therefore
employ a modification of the Stretch Lemma: in our simulations below we will stretch the input to
occupy every {\it fourth} cell, using blocks of four cells. The first three cells in such a block 
represent one cell each from the input, scratch and output tapes, while the fourth simply 
holds a flag that identifies the left-most block of
cells, and it is never written on during the computation or changed in any way. By this means, all the information will be available in order for us to correctly simulate the operation of any three-tape machine computation on a one-tape machine.

\proofof \Almost: The reverse direction clearly holds, so we prove the 
forward direction. Suppose $f$ is computable by the program $p$. We
will design an algorithm for computing $1\concat f$ on a one-tape
machine. The algorithm, following the three phases of our general
strategy above, will stretch the input, then simulate $p$, and finally
compress the output before halting. 

Thus, from the {\it start} state our
algorithm takes a real input $a$ and stretches it in just over
$\omega^{2}$ steps as in the Stretch Lemma \Stretch, but with the
input occupying every fourth cell, and with two
additional cells at the start of the tape to use as flags for keeping
track of which phase of computation we are in.   
After $\omega^{2}$ many steps, the tape will read
$\<00a_{0}101a_1000a_2000\cdots>$.  Recall that the $1$ in the fourth
cell is used in the algorithm of the Stretching Lemma to signal that
the stretching is done, and the $1$ in the sixth cell is the fourth cell in the first block of 
four cells, marking that that block of four is the first block. The fourth cell of 
every subsequent block of cells remains $0$ throughout the computation.

We now erase the $1$ in the fourth cell, and write a $1$ in the
second cell, to indicate that we have entered phase two of the
simulation.  Thus, a few steps
after $\omega^2$, the tape reads $\<01a_0001a_1000a_2000\cdots>$. 
The $1$ remains in the second cell to remind the machine at each limit
stage that it is in phase two of the simulation, the phase in which it
is directly simulating the computation of $p$ as in the Simulation
Lemma \Simulate, but without modifying the first two cells.  

When this simulation is
complete (that is, when the simulated program reaches its halt state)
our algorithm erases the $1$ in the second cell and writes a $1$ in
the first cell, signalling that phase three has begun. At this point, the tape 
reads $\<10x_0y_0z_01x_1y_1z_10\cdots>$ where $f(a)=z$ and $x$ 
and $y$ are the contents of the simulated input and scratch tapes at
the end of the computation of $f(a)$. The machine now compresses the
output in $\omega$ steps, so that by the next limit stage the tape
reads $\<1z_0z_{1}z_{2}\cdots>$.  Then, the machine notes the
$1$ in the first cell and, recognizing that this is the first 
limit stage at which it has seen a $1$ in that cell, knows that it
has completed the simulation. It is therefore able to halt with output $1\concat
z$. Since $z=f(a)$, this is precisely the output for $1\concat f(a)$,
as desired.\qed

With the idea of this theorem we can now prove the one-tape computability of many
different functions. 

\theorem. If the range of $f\from\R\to\R$ is not dense, then $f$ is
computable if and only if it is one-tape computable. 
\ref\Notdense

\proof If the range of $f$ is not dense, then it omits an open set,
and so there is a finite string $\sigma$ which is not an initial
segment of $f(a)$ for any $a$. Consider the algorithm which
computes $1\concat f$ on one tape according to the procedure described
in \Almost, except that $\s$ is appended to the initial part of the
tape all throughout the computation (with the algorithm checking at
every limit to make sure $\s$ is still there). Thus, the tape reads
$\s\concat 1\concat z$ at the conclusion of this procedure. At this point, for
the first time, the algorithm erases the initial $\s\concat 1$ and
shifts the true output $z$ over by $|\s\concat 1|$ many steps, leading
to a tape with just $z$ on it. Since $z=f(a)$ does not have $\s$ as an
initial segment, this is the first limit stage at which the first few
cells on the tape differ from $\s$, and the algorithm can recognize
this. Knowing now that it has completed the computation, the
algorithm halts.\qed 

The next two results follow immediately.

\corollary.  A function $f\from\R\to \N$ is computable if
and only if it is one-tape computable.

\corollary.  A function $f\from\R\longrightarrow\{0, 1\}$ is
computable if and only if it is one-tape computable. Consequently, a
set of reals is decidable if and only if it is one-tape decidable. 
\ref\Characteristic

The following corollary to Theorem \Almost\ demonstrates that the
condition given in Theorem \Notdense\ --- that the range of the 
function is not dense --- is not necessary for the function to be
one-tape computable.

\corollary.  A function $f\from\R\longrightarrow\Q$ is computable if
and only if it is one-tape computable.

\proof Every element of $\Q$ can be coded in a canonical way with an
element of $\N$, by means of some appropriate G\"odel coding, and so
to compute a function $f\from\R\to\Q$ one first computes the analogous
function into $\N$ and then transforms the element of $\N$ in finitely
many steps into the corresponding element of $\Q$.\qed 

We conclude this section with a theorem showing that two models of
infinite time Turing machines give rise to the same degree
structure. We assume that the one-tape machines are augmented in the
same way as the three-tape machines with oracle tapes. 

\theorem. One-tape infinite time Turing machines give rise to the same
notions of relative computability $A\ilte B$ and computable
equivalence $A\iequiv B$ as do the regular infinite time Turing
machines. 

\proof This theorem follows by simply relativizing Corollary
\Characteristic\ to oracles. To do so, one relativizes the Simulation
Lemma \Simulate\ to oracles. In the case of real oracle, where the
oracle is written out on the oracle tape, one organizes the one
computation tape in blocks of four rather than three, so that each
block has cells for the oracle, input, scratch and output
tapes. During the set-up phase, one must stretch both the input and
the oracle to occupy the appropriate cells in the simulation, and having thus 
copied the oracle tape into the simulation form, one thereafter ignores the 
actual oracle tape. In the
case of an oracle which is a set of reals, one has a blank oracle tape
on which reals can be written and queries made about their membership
in the oracle. To simulate this on a one-tape machine with an oracle
tape, one divides the computation tape into blocks of five, representing 
cells for the input, scratch, output and oracle tapes as well as one additional 
information cell. Whenever a query is made in the simulation, the algorithm copies
the contents of the simulated oracle tape to the actual oracle tape to
make the query, using the fifth cell to keep track of where the head was 
before this and also to signal at a limit that the algorithm has just performed this 
operation. In summary, either kind of oracle can be simulated with a
one-tape machine augmented with an oracle tape. To prove the theorem
at hand, now, recall that the relation $A\ilte B$ holds when the
characteristic function of $A$ is computable with oracle $B$. Since
the characteristic function of a set has range in $\set{0,1}$, this
function will be one-tape computable from $B$ by the relativized
argument of \Characteristic. And once the notions of $\ilte$ are seen
to be the same for the two kinds of machines, it follows that the
notions of $\iequiv$ are also the same.\qed

\section The one-tape machines are not fully powerful

After the results of the previous section, we appear to be on the verge of
showing that every computable function $f$ is one-tape computable. If
$f$ is computable, then we have proved that $1\concat f$ is one-tape
computable; all that remains is to remove the extra $1$ from the
front and, in $\omega$ many steps, shift the rest of the cells to the
left by one.  The problem is how, after shifting the cells of $1\concat f(a)$ to
the left to form $f(a)$ by the next limit stage, would we know to
halt?  At that stage, the tape contains $f(a)$, essentially an
arbitrary real.  There can be no flags to signal that 
we're done shifting.  In fact, this output real could have appeared on
the tape at some earlier limit stage, and the algorithm be caught in a
loop. Behind this problem lies a surprising result. 

\theorem. There is a computable function which is not one-tape computable.
\ref\Counterexample

\proof
We will construct the desired function by diagonalizing against all
one-tape machines. Before doing so, let us introduce some helpful
terminology. If a program halts after a limit stage, then it does so
because a finite initial segment of the tape supports a halting
computation from the limit stage. So let us say that a finite sequence
$\sigma\in 2^{<\omega}$ of length $|\sigma|$ is a {\df halting string}
for a one-tape program $p$ if the program $p$, encountering $\sigma$
on the tape at a limit stage with its head on the first cell and in
the limit state, halts in less than or equal to $|\sigma|$ steps. We
refer to the corresponding finite output $\tau$ with
$|\tau|=|\sigma|$, the result of the computation, as the corresponding
{\df halting string output} for $p$. If $x$ is in the range of the
function computed by $p$, and that computation halted after a limit
stage, then all sufficiently long initial segments of $x$ are halting
string outputs. (But if a program halts in finitely many steps, that
is, before reaching any limit stage, then the output may have nothing
to do with halting string outputs.) It is easy to see that any
extension of a halting string is a halting string, and any extension
of a halting string output is a halting string output.  Thus, the set of
reals that do not extend a halting string output for $p$ is closed,
and, if non-empty, contains a lexically least element $x_p$ (the
left-most branch through the tree of all non-halting string outputs).

Recall from \cite[HamLew$\infty$a] that a {\df writable} real is one which is the output of a supertask program on input $0$. By the results of \cite[HamLew$\infty$a], we may fix a writable real $u$ which is not writable by any machine in fewer than
$\omega+\omega$ steps.  We may also fix a computable enumeration of all one-tape
Turing machine programs $p_{0}, p_{1}, p_{2}, \ldots$. Define now the partial function $f$ by $f(0)=u$ and, for $n>0$, $f(n+1)=x_{p_n}$, if it exists. 
The function is not defined for non-integer inputs.

We claim first that $f$ is computable by an ordinary three-tape
program.  On input $0$, we instruct the machine to write $u$.  On
input $n+1$, we instruct the machine to construct the list of all
halting outputs of the program $p_{n}$, and then output the lexically
least branch through the tree which is the complement of that set, namely, $x_{p_n}$, if it exists. Otherwise, our algorithm does not halt. 

We conclude the proof by showing that $f$ is not one-tape computable.
Assume towards a contradiction that it is, by some program
$p_{n}$. Let
$x$ be the real written on the tape at stage $\omega$ by $p_{n}$ on input
$0$.  By assumption, $p_{n}$ computes $f(0)=u$ correctly, and by the
choice of $u$ this must take at least $\omega+\omega$ many steps. Thus, no
initial segment of $x$ can be a halting string, for otherwise the
computation would halt by some stage $\omega+k$. In particular, there
are infinitely many strings which are not halting strings for
$p_n$. And since each halting string of length $n$ gives rise to
exactly one halting string output of length $n$, it follows also that
there are infinitely many strings which are not halting string outputs
for $p_n$. Thus, by K\"onig's Lemma, the tree of all such strings has
a branch, and so it has a lexically least branch. That is, $x_{p_n}$
exists. Consequently, by definition, $f(n+1)=x_{p_n}$.  

Since we have assumed that $p_n$ computes $f$, it must be that on
input $n$ the program $p_n$ gives output $x_{p_n}$. But no initial
segment of $x_{p_n}$ is a halting output string for $p_n$, so by the
remarks in the first paragraph of this proof, this computation cannot
be the result of an infinite computation. So the program must have
halted in some finite number of steps $k$. But in a finite
computation, the machine only has a chance to view the first $k$
digits of the input before it halts. By manipulating the input $n+1$
past its $k^\th$ digit, we can create a new input $z$ which is not an
integer but which leads to the same halting computation with program
$p_n$. This contradicts the fact that the domain of $f$ is contained
in the integers. We conclude that $f$ is not computable by a one-tape
machine.\qed 

\corollary.  The class of one-tape computable functions is not closed
under composition.  The closure of the class of one-tape computable
functions under composition is exactly the class of all infinite time
computable functions. 

\proof The proof of the Overly Hopeful Theorem shows that every
computable function can be expressed as a composition of one-tape
computable functions.  Since not every function is one-tape
computable, the class of one-tape computable functions is not closed
under composition.\qed 

The idiosyncratic nature of infinite time Turing machines with only
one tape is the simple result, we believe, of cramped working
space.  The situation is like that faced by a great artist painting a masterpiece on 
a vast cathedral floor;  before completing the work, the artist finds himself
with nowhere to stand. Similarly, the one-tape machines find that if
the whole tape is to become the output of the computation, there is no
room left for flags to signal side information about the computation,
such as when it is complete. 

Carrying this idea further, let us augment the one-tape machines with
a scratch pad consisting of one cell. The machine can
read from and write on this scratch pad cell, but the value of the
scratch pad does not become part of the output. The scratch pad 
provides our artist with a place to stand in the end, so to speak, without
disturbing the masterpiece of the output itself. 
$$\hbox{\vtop{\baselineskip=0pt\parskip=0pt\cell0\vskip1pt\hbox{\sixrm\hskip.5pt scratch}\hbox{\sixrm\ \ pad}}}\hskip 5pt\vtop{\baselineskip=0pt\parskip=0pt\hbox{\cell1\cell0\cell1\lower2pt\vbox{\hbox{\strut\ninerm\hskip3pt Head}\boxit{\cell1}}\cell0\raise5pt\hbox{\ $\cdots$}}\vskip2.5pt\hbox{\sixrm\quad input/output tape}}$$
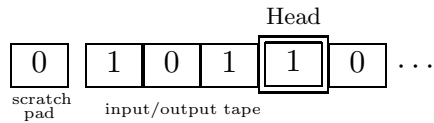
$$\hbox{\eightrm A one-tape infinite time Turing machine with a one-cell scratch pad}$$

\theorem. One-tape infinite time Turing machines with a one-cell scratch pad are fully as powerful as the three-tape machines. 
\ref\ScratchPad

\proof A machine with a one-cell scratch pad works in effect just like
a one-tape infinite time Turing machine, except that the scratch pad cell
value is a part of neither the input nor the output. Given now a
computable function $f$, we may compute $1\concat f$ on this machine,
treating the scratch cell as if it were the first cell on an ordinary one-tape machine. 
At the end of the computation we have
$1\concat f(a)$ on the tape. Since the initial $1$ sits on the scratch
pad cell, the actual output is $f(a)$, as desired.\qed 

We can similarly modify the notion of infinite time Turing machines to
include $n$ tapes for any $n$, and obtain the following corollary. 

\corollary. For any $n>1$, a function is $n$-tape computable if and
only it is computable.

\proof The regular infinite time Turing machines, with $n=3$, can
simulate the $n$-tape machines by means of G\"odel
pairing. Conversely, the previous theorem shows that having a one-cell
scratch pad is sufficient to compute all infinite time computable
functions.\qed 

We take the results of this section --- most notably, that the class of one-tape computable functions 
is not closed under composition --- to show that the model of one-tape infinite time Turing machines 
is not the right notion, and does not provide the right model of supertask 
computation. But when one augments such a machine with a scratch tape of any size, even a scratch pad consisting of 
just one cell, the results show that one arrives at the same robust notion of computation as the original definition of the infinite time computability. Indeed, we take this as an affirmation of the robustness of the orginal multi-tape definition. Perhaps our results show 
that the simplest kind of machine leading to the same full concept of supertask computation is an infinite time Turing machine with one tape augmented by a one-cell scratch pad. And we propose, especially for those wanting to work with a one-tape model, that these 
one-tape infinite time Turing machines with scratch pad provide the correct one. 

\section Clockable ordinals

In this final section we investigate the relationship between the
clockable ordinals and their one-tape clockable counterparts. Recall
from \cite[HamLew$\infty$a] that an ordinal $\alpha$ is {\df
clockable} when there is an infinite time Turing machine program which
halts on input $0$ in exactly $\alpha$ many steps. Similarly, an
ordinal is {\df one-tape clockable} when it is the length of a
one-tape infinite time Turing machine computation. Results in
\cite[HamLew$\infty$a] include, for example, the facts that every
clockable ordinal is countable, that there are gaps in the clockable
ordinals of length unbounded in $\gamma$, the supremum of the
clockable ordinals, and that the first gap begins at 
$\omega_1^{\sixrm CK}$ and has length $\omega$. In fact, in 
\cite[HamLew$\infty$a]
it is proved that no admissible ordinal is clockable. Philip Welch
proved in \cite[Wel$\infty$a] that $\gamma$ is also equal to the
supremum of the writable ordinals, the ordinals coded by a writable
real, and consequently every clockable ordinal is writable. Let us say
that a clockable ordinal $\alpha$ {\df ends a gap} in the clockable
ordinals when there is an interval $[\beta,\alpha)$ containing no
clockable ordinals; the least such $\beta$ is the corresponding {\df
gap-starting} ordinal. 

Since the ordinary three-tape infinite time Turing machines can
directly simulate, step-for-step, the machines with only one tape, it
follows that every one-tape clockable ordinal is clockable. The
question is whether the converse holds.  

\theorem.  Every clockable ordinal that is not one-tape clockable ends
a gap in the clockable ordinals. 
\ref\Clockable

The theorem is an immediate consequence of the next two lemmas.

\quiet\lemma. Every clockable successor ordinal is one-tape
clockable. 

\proof  Suppose $\a+n$ is a clockable successor ordinal, where $\a$ is
a limit ordinal and $n\geq 1$ is a positive integer. It suffices for
us to show that $\a+1$ is one-tape clockable, since any computation
can be prolonged finitely many steps by means of counting through
extra states. By the Speed-Up Lemma of \cite[HamLew$\infty$a], we know
that $\alpha$ itself is clockable. We will clock $\a+1$ with a
one-tape machine by simulating the computation of a three-tape machine
clocking $\a$. 

Fix a program $p$ clocking $\alpha$ and let $i, j, k\in\{0,1\}$ be the
digits appearing in the first cell of the input, scratch and output tapes
at stage $\alpha$ in the computation of $p$ on input $0$. Since this
computation halts at stage $\a$, it must be that $\alpha$ is the first
limit stage at which $i$, $j$ and $k$ appear in those cells.  We will
simulate $p$ with a 
one-tape program $q$ by anticipating the appearance of $i$, $j$ and
$k$ in essentially the same manner as the proof of the Speed-up Lemma
in \cite[HamLew$\infty$a], through the use of flags located on the
first two cells of the tape. Let us refer to the first of these flag
cells as the $0$-flag and the second as the $1$-flag.  

We simulate $p$
with a one-tape program $q$ as in the Simulation Lemma \Simulate,
leaving room in the front of the tape for the two flag cells and
whatever additional space is necessary for book-keeping. Since the
input is $0$, there is no need for the stretching phase. After each
step of the $p$ computation simulation, we return the head to 
the first three simulated cells and compare their contents with $\<
i, j, k >$.  If each of these three simulated cells contains a $0$
when the corresponding digit in $\< i, j, k >$ is $0$, then we
set the $0$-flag to $0$.  Otherwise, we set the $0$-flag to $1$.  We
then examine 
those of the three cells for which the corresponding digit in $\< i,
j, k>$ is $1$.  We flash the $1$-flag if each of these cells
contains a $1$ or has displayed a $1$ at some point since the last
flash.  The point of this procedure, which is easy to verify, is that
at a limit stage the $0$-flag is $0$ and the $1$-flag is $1$ precisely
when the 
first three simulated cells agree with $\< i, j, k>$. And since the
simulation procedure simulates $\omega$ many steps of the computation
of $p$ in $\omega$ many steps, catching up as it were at every limit
stage, this means that the two flags are set to $0$ and $1$,
respectively, for the first time at a limit stage, exactly at
$\alpha$. So, by checking the status of these two flags at every limit
stage, the algorithm will be able to halt precisely at $\alpha+1$.\qed

\quiet\lemma. Every clockable limit of clockable ordinals is one-tape
clockable.

\proof  Suppose that $\alpha$ is a clockable limit of clockable
ordinals. We will modify the algorithm of the previous lemma in such a
way as to avoid the need for the $0$-flag, and therefore the need to
take the extra step checking that flag at $\alpha$. The modified
algorithm will therefore be able to halt in exactly $\alpha$ many
steps. 

Let $\< i, j, k>$ be, as above, the contents of the first cells on
each of the three tapes in the computation clocking $\alpha$. If one
or more of these cells displays a $0$ at $\alpha$, then it must have
been $0$ from some point on before $\alpha$.  Since $\alpha$ is a
limit of clockable ordinals there is some clockable $\beta<\alpha$ by
which those cells have stabilized to $0$, being $0$ for the duration
from $\beta$ up to $\alpha$. We now simultaneously simulate the
program $p$ and a program clocking $\beta$ on a one-tape machine.  We
leave the
first cell of the tape available for use as the $1$-flag, but we begin
flashing it only after stage $\beta$ has been reached. Thus,
the $1$-flag will contain a $1$ at a limit stage for the first time at
stage $\alpha$. By checking this flag at every limit stage, the
machine will be able to halt in exactly $\alpha$ many steps.\qed 

\theorem.  There exist clockable ordinals which are not one-tape
clockable.   In particular, any ordinal which ends a gap of compound
limit length is not one-tape clockable.  
\ref\Notonetape

\proof  We prove the second sentence first, that any ordinal which
ends a gap of compound limit length is not one-tape clockable. Suppose
to the contrary that $\alpha$ ends a gap of compound limit length and
is one-tape clockable. Since the length of the gap leading up to
$\alpha$ is a limit of limit ordinals, it follows that the ordinal
$\alpha$ itself is a limit of limit ordinals. Consider now the
one-tape computation which clocks $\alpha$. Since $\alpha$ is a limit
ordinal, we know that at stage $\alpha$ in this computation, just
before halting, the head is on the first cell of the tape, the machine
is in the {\it limit} state, and there is either a $0$ or a $1$ on the
first cell of the tape. The cell actually cannot display a $0$ at
stage $\alpha$, because then the cell would have been $0$ from some
point on before $\alpha$, and since $\alpha$ is a limit of limit
ordinals, it would have been $0$ at some limit stage before $\alpha$,
therefore causing the computation to halt at that earlier time. Thus,
at stage $\alpha$, the first cell must be $1$. It follows that the
cell must have displayed a $1$ unboundedly often in $\alpha$ and in
fact $\alpha$ must be the $\omega^\th$ time that this cell displays
$1$, since otherwise the computation would have halted earlier, at the
$\omega^\th$ instance. It is easy now to see that $\alpha_n$, the
$n^\th$ time this first cell is $1$ during the computation, is
clockable. Furthermore, since the $\alpha_n$ are unbounded in
$\alpha$, it must be that $\alpha$ is a limit of clockable ordinals,
contradicting our assumption that it ends a gap. Thus, the second
sentence of the theorem is proved; any ordinal which ends a gap of
compound limit length is not one-tape clockable. 

We now prove the first sentence by showing that such ordinals
exist. We could simply show that the first gap of size at least
$\omega^2$ has size exactly $\omega^2$, in order to conclude that gaps
of compound limit ordinal length exist.  It is not much more
difficult, however, to prove that if $\d$ and $\b$ are clockable
ordinals, with $\b$ a limit, the first gap above $\d$ of 
size at least $\b$ has size exactly $\b$. When $\b$ is a compound
limit ordinal, the theorem follows. Let $\b'$ be such that
$\omega+\b'=\b$.  If $\omega^2\leq\b$, then of course $\b'=\b$.  In
any event, $\b'$ is clockable.  

Fix $\d$ and $\b$ and observe that by the results of
\cite[HamLew$\infty$a] there is a first gap above $\d$ of size 
at least $\b$.  We will design an algorithm to recognize this
gap, and use this algorithm to show that the gap cannot be 
longer than $\b$. In the manner of many of the arguments of
\cite[HamLew$\infty$a], by means of G\"odel-coding we imagine that the
scratch tape of a three-tape supertask machine is divided into
$\omega$ many scratch tapes, each used to simulate the operation of
one of the infinitely many infinite time Turing machine programs on
input $0$, keeping careful track of which programs have halted. The
simulation procedure is set up so that $\omega$ many steps of
simulated computation are carried out for each program in $\omega$
many steps, so that the process catches up at every limit stage. We
also reserve room for two additional clocks, one counting to $\d$ and
one counting to $\b'$.  

While simulating the computation of all programs, we count to $\d$.  After $\d$,
every time we find that one of the simulated programs halts, we flash 
a flag on the first cell of the tape.  Any limit stage, therefore, at
which this flag is $1$ must be a limit of clockable ordinals (and at
this point we reset the flag to $0$). If at a limit the flag is $0$,
that means that the simulated computations are not halting, and we are
in or have just finished a gap.  Observe that this flag is first $0$
at a limit at $\omega$ many steps past the beginning of the gap,
because the beginning of the gap is, of course, a limit of clockable
ordinals, and the algorithm must check all the programs before being
sure that none have halted at that stage. The algorithm in effect
recognizes gaps $\omega$ steps past their occurence.  

At any limit stage at which the flag is $0$, that is, while we are
inside a gap, we use the $\b'$ clock counting to $\b'$ --- while
continuing to search for halting programs --- in order to determine
the length of the gap. We reset this clock if the gap runs out before the 
$\beta'$ clock, and at limits of such resettings. If it happens that the $\b'$ clock finishes
before we find another halting program, then, because of the initial
lag of $\omega$ steps in recognizing the gap, the gap has size at
least $\omega+\b'=\b$, so we halt.

Since we discover this gap of size at least $\b$ within
finitely many steps of the $\b^\th$ ordinal past the start of the gap,
this means that our algorithm will halt within finitely many steps of
the $\b^\th$ ordinal past the start of the gap.  Consequently, that
gap has size less than $\b+\omega$, and so by the Speed-up Lemma of
\cite[HamLew$\infty$a], it must have size exactly $\b$.\qed 

The only ordinals remaining whose one-tape clockability status is in
question are the clockable ordinals that end gaps of simple limit
length, that is, that end gaps of length $\b+\omega$ for some
$\b$. The theorem below, which generalizes readily, shows that many of
these are one-tape clockable. 

\theorem.  If $\a$ is the least gap-ending ordinal above a given clockable
ordinal, then $\a$ is one-tape clockable. More generally, if $\a$ is the
least ordinal above $\delta$ which ends a gap of length at least
$\b+\omega$ where $\delta$ and $\b$ are clockable, then $\a$ is one-tape
clockable. Indeed, $\b$ need only be one-tape writable in time before
the gap in question. 

\proof  Suppose $\a$ is the least gap-ending ordinal above the
clockable ordinal $\delta$.  By a result of \cite[HamLew$\infty$a], $\a$
ends a gap of length $\omega$.  Thus, the first sentence follows from
the second by letting $\beta$ be finite.  We prove this case first.

We describe an algorithm with a single master
flag which signals when the machine should halt at a limit stage. 
Consider the algorithm from Theorem~\Notonetape\ which simulates all three-tape 
infinite time Turing machine computations on input $0$, keeping
careful track of which programs halt, while simultaneously 
counting to $\delta$. We keep the master flag set to $1$ until stage
$\delta$, but after stage $\delta$ we set it to $0$ after every limit stage,
and flash it each time we find that a simulated program has
halted. The first time this master flag is $0$ at a limit stage will
therefore be $\alpha$, since that is the end of the first gap after
$\delta$.  By placing this flag on the first cell of the tape and
checking it at every limit stage, we can halt right at
stage $\alpha$, as desired.  

We now prove the second sentence in the case $\beta$ is infinite by 
designing a program that causes a one-tape machine to halt at the
first ordinal $\alpha$ such that:
\points 1. $\alpha>\delta$.\cr
	2. $\alpha$ lies in or ends a gap (in fact, it will end the gap).\cr
	3. $\alpha$ is the $(\beta+\omega)^\th$ ordinal past the start
of this gap.\cr   

\noindent
Our algorithm has a single master flag which will signal when all
three conditions are met, since this $\alpha$ must end the first gap above 
$\delta$ of size at least $\beta+\omega$, the theorem will be proved.

The only awkward point is noticing when we have satisfied condition
$(3)$. Searching for gaps in the fashion of the previous
Theorem~\Notonetape, we 
will not recognize that we have found a gap until we have reached  
the $\omega^\th$ ordinal after its start. In order to count to the
$(\beta+\omega)^\th$ ordinal past the start of the gap, we will fix
$\beta'$ such that $\beta=\omega+\beta'$, and count to
$\beta'+\omega$.  Since $\beta$ is clockable, so is $\beta'$.    

We now describe the algorithm in detail.  From the {\it start} state,
the machine places a $1$ in the first three cells.  The first cell is
our master flag; we halt when it is $0$ at a limit stage.  The second
cell records when we have reached $\delta$, and the third cell records
when we have found a gap.

While simulating all computations on input $0$, the machine begins by also 
counting to $\delta$.  After reaching $\delta$, it enters a $0$ in the
second cell to remind us at each limit stage that we have passed $\delta$.  
The machine now continues to simulate all three-tape infinite time Turing machine
computations on input $0$, keeping careful track of which 
programs halt. We write a $0$ in the third cell at each limit stage, and
flash the third cell each time we find that a simulated program has
halted.  Any limit stage at which the third cell contains a $1$ is a
limit of clockable ordinals.  If at a limit this flag is $0$, then the
simulated computations are not halting, and we are in or have just
ended a gap.    

Observe that the first limit at which the second and third cells are
both $0$ occurs $\omega$ steps past the beginning of the first gap
above $\delta$.  At this point, the machine begins to count to $\beta'$
while continuing to check for halting programs.  
If the gap ends before the $\beta'$ clock runs out, then it was too short, 
and the machine resets
the $\beta'$ clock and goes on searching for other gaps. The $\beta'$ clock 
is also reset at limits of such resettings. Otherwise, after
the $\beta'$ clock runs out, the machine sets the master flag to $0$. 
While the master flag is $0$, if we find another program to halt, then we know that the 
gap was too short (having length only $\b$), and we reset the master flag and the third cell to $0$,
resets the $\beta'$ clock, and continue to search for gaps. Otherwise, if
no simulated program halts, then the machine will halt at the next limit
stage, which is $\alpha$. This concludes the proof of the second sentence.  

We now prove the final remark --- the case that $\b$ is not
necessarily one-tape clockable, but is one-tape writable in time
before the gap in question.  That is, we assume that before the
start of the gap that $\alpha$ ends, we can write a real
coding a relation on $\omega$ with order-type $\b$.  This real can be
used as a clock for counting to $\b$ by gradually erasing its initial
segments. One tick of the  clock consists of finding and marking as
{\it deleted} the least element from the field of the relation. One
can execute $\omega$ many ticks of the clock in $\omega$ many steps by
finding and marking as deleted the $\omega$ least members of the field
of the relation. Furthermore, one can tell at a limit stage that the
real has already been completely deleted with a single master flag by
flashing such a flag each time the least element of the relation (in
the natural number order of the G\"odel codes) is deleted. This flag
will be on for the first time at a limit stage when the relation has
been entirely deleted. Thus, by modifying the previous algorithms to
use this real as a clock, one obtains the result in the case that $\b$
is writable in time before $\a$, as desired.\qed  

So we know that many of the gap-ending ordinals that end gaps of
simple limit length are one-tape clockable, and we know of no such
ordinals that are not one-tape clockable. The following questions
remain open.  

\question Open Question. Exactly which clockable ordinals are not
one-tape clockable? In particular, are the clockable non-one-tape
clockable ordinals exactly the ordinals which end gaps of compound
limit ordinal length?

There is a one-tape model of computation having both the same notion of computability and the same clockable ordinals as the ordinary infinite time Turing machines. Specifically, in addition to augmenting a one-tape machine with a one-cell scratch pad, we propose to use a double-sized head, capable of reading two cells at once. $$\hbox{\vtop{\baselineskip=0pt\parskip=0pt\cell0\vskip1pt\hbox{\sixrm\hskip.5pt scratch}\hbox{\sixrm\ \ pad}}}\hskip 5pt\cell1\cell1\cell0\lower2pt\vbox{\hbox{\strut\hskip5pt\sixrm Double Head}\boxit{\cell1\cell0}}\cell0\cell0\cell1\raise5pt\hbox{\ $\cdots$}$$
Since the input and output fill the entire tape, except for the one-cell scratch pad, such a model is not automatically the same as a two-tape machine. Nevertheless, we have the following theorem:

\theorem. One-tape double-head infinite time Turing machines with scratchpad lead to the same class of computable functions and the same set of clockable ordinals as the ordinary infinite time Turing machines. 

\proof The scratch pad is sufficient to compute the same computable functions by Theorem \ScratchPad. The double-head
allows the machine to view both the $0$-flag and the $1$-flag in the argument of the first lemma of \Clockable, thereby allowing the machine to 
halt right at $\alpha$ in that argument, and so these machines will have the same clockable ordinals.\qed

We would like now to close this paper by turning to the question of the efficiency of
one-tape supertask machines.  Corollary \Characteristic, asserting
that every infinite time decidable set is one-tape decidable, relies
on our three-step algorithm for simulating a three-tape computation
with a one-tape machine. By analyzing the time each of these steps
takes, we obtain the following result.

\theorem Efficiency Theorem. For every infinite time Turing machine
program $p$ deciding membership in a set $A\of\R$, there is a one-tape
program $q$ deciding $A$ such that if $p$ takes $\a$ many steps to
decide whether $a\in A$, then $q$ takes $\omega^2+\a+\omega$ steps to do so. 
If $p$ is sufficiently tidy, then the computation of $q$ can
be arranged to take only $\omega^2+\a+1$ many steps. (note that for
$\a\geq\omega^3$, this is the same as $\a+1$)
\ref\Efficiency 

\proof On input $a$, our strategy first called for stretching the
input in $\omega^2+1$ steps. Then, the algorithm simulates the
operation of $p$ with seven-steps-for-one, and then compresses the
output in $\omega$ additional steps. This takes
$\omega^2+1+7\cdot\a+\omega=\omega^2+\a+\omega$ many steps in all. If
the program $p$ is tidy in the sense that it leaves nothing except
$0$s on the input and scratch tapes at the end of the computation,
then the simulation can be performed in $\omega^2+\alpha+1$ many
steps, because the output is $0$ or $1$, and so with a tidy
computation the final compression phase can be omitted. To get the
$+1$ in this case, one should put the output cell on the first cell of
the tape, with the master halt flag on the second cell (rather than on
the first cell as in Theorem \Almost), so that the program needs only
one additional step after $\a$ to check it and halt.\qed

One naturally wonders whether this bound can be improved. Certainly
one cannot expect in general to decide $A$ in $\a$ many steps, because
$\a$ may be clockable but not one-tape clockable. Because of this,
$\a+1$ seems the best possible general bound. 

\question Open Question. Can every supertask computation in $\a$ steps ($\a$ infinite) be uniformly simulated on a one-tape machine in $\a+1$ many steps?

One might hope to answer this question by improving the $\omega^2$ term in our Theorem \Efficiency, the term which arises from applying the stretch function to the input. 

\question Open Question. How long does it take to compute the stretch function with a
one-tape machine?  

A partial answer to the question above is given in the following
theorem.  

\theorem.  The stretch function is not computable in $\omega$ steps by a
one-tape machine, even by one augmented with a scratch pad of any finite size. 

\proof The stretch function, defined by 
$s:\<x_0x_{1}x_{2}x_{3}\cdots>\mapsto\< x_000x_{1}00x_200\cdots>$,
stretches the input to occupy every third cell. For convenience, we
refer to the other cells 
as the $0$-cells.  Assume towards a contradiction that $s$ is computable in $\omega$ many steps on every input by a
one-tape machine augmented with a scratch pad of $p$ many cells. Choose $m$ large enough 
so that $m^22^p<2^m$ (for example, if $p=0$, then it suffices to take $m>4$). Let us also assume that $m$ is at least as large as the number
of states in the program. Since $\R=2^\w$ is compact, there is a
sufficiently large $n$ such that for each real input, there is some
$k\leq n$ such that at the $k^\th$
stage of computation, each of the first $m$ many $0$-cells has a $0$
written in it. Associate to each input of length $n$ the machine's
configuration --- cell contents, state, and head position --- at the
first 
such stage.  We now count configurations. There are at most $m$ many
states, $m$ 
many possible head positions (since a $0$ has just been written on one of
$m$ many cells), at most $2^p$ possible strings filling the scratch pad and at most $2^{(n-m)}$ possible strings filling the
remaining cells. Thus, each of the $2^n$ input strings of length $n$
is associated to one of $m^{2}2^{n-m}2^p$ many configurations. Since this is less than $2^n$, we
conclude that there are two distinct strings $t$ and $u$ of length $n$
leading to the same configuration.  It now follows that for any real
$y$, the inputs $t\concat y$ and $u\concat y$ will lead eventually to
the same computation. This contradicts the fact that $s$ is a
one-to-one function.\qed  

The argument does not seem to generalize easily to $\omega+\omega$.

\section Bibliography

\nopagenumbers
\parindent=0pt
\newbox\Article
\newbox\Journal
\newbox\Author
\newbox\Vol
\newbox\No
\newbox\Year
\newbox\Page
\newbox\Book
\newbox\Publisher
\newbox\Pubaddr
\newbox\Key
\newbox\Editor
\newbox\Comment
\newbox\Note
\def\entry#1#2\par{\item{#1\quad}\hskip-1.1em#2\par}
\def\article#1{\setbox\Article=\hbox{\sl #1, }}
\def\journal#1{\setbox\Journal=\hbox{\rm #1 }}
\def\author#1{\setbox\Author=\hbox{\sc #1, }}
\def\vol#1{\setbox\Vol=\hbox{\bf #1 }}
\def\no#1{\setbox\No=\hbox{no. #1 }}
\def\year#1{\setbox\Year=\hbox{\rm({\oldstyle #1}) }}
\def\page#1{\setbox\Page=\hbox{\rm p. #1 }}
\def\book#1{\setbox\Book=\hbox{\it #1, }}
\def\publisher#1{\setbox\Publisher=\hbox{\rm #1, }}
\def\pubaddr#1{\setbox\Pubaddr=\hbox{\rm #1, }}
\def\key#1{\setbox\Key=\hbox{#1}}
\def\editor#1{\setbox\Editor=\hbox{\rm(#1, Ed.) }}
\def\comment#1{\setbox\Comment=\hbox{\rm #1}}
\def\note#1{\setbox\Note=\hbox{\rm #1 }}
\def\ref#1\par{\smallskip{#1
  \entry{\ifhbox\Key\unhbox\Key\else[\ ]\fi}%
  \unhbox\Author\unhbox\Note
  \ifhbox\Book \unhbox\Book\unhbox\Publisher\unhbox\Pubaddr
               \unhbox\Editor\unhbox\Page\unhbox\Year\unhbox\Comment
  \else \unhbox\Article\unhbox\Journal\unhbox\Vol\unhbox\No\unhbox\Editor
        \unhbox\Page\unhbox\Year\unhbox\Comment\fi\par}}

\tenpoint\tightlineskip\leftskip=.5in

\ref
\author{Joel David Hamkins and Andrew Lewis}
\article{Infinite Time Turing Machines}
\journal{to appear in the Journal of Symbolic Logic}
\key{[HamLew$\infty$a]}

\ref
\author{Joel David Hamkins and Andrew Lewis}
\article{Post's problem and infinite time Turing machines}
\journal{submitted to the Archive for Mathematical Logic}
\key{[HamLew$\infty$b]}

\ref
\author{Robert Soare}
\book{Recursively Enumerable Sets and Degrees}
\publisher{Springer-Verlag}
\year{1980}
\key{[Soa]}

\ref
\author{Philip Welch}
\article{The lengths of infinite time Turing machine computations}
\journal{to appear in the Bulletin and Journal of the London Math. Soc.}
\key{[Wel$\infty$a]}

\ref
\author{Philip Welch}
\article{Eventually infinite time Turing degrees: infinite time decidable reals}
\journal{to appear in the Journal of Symbolic Logic}
\key{[Wel$\infty$b]}

\ref
\author{Philip Welch}
\article{Friedman's trick: minimality in the infinite time Turing degrees}
\journal{submitted to ``Sets and Proofs', vol. 2, 1997 Proc. European
Meeting of the ASL, Leeds, (CVP LMS Lecture Notes in Maths Series 1998)}
\key{[Wel98]}

\bye